\documentclass[twocolumn,11pt]{article}
\usepackage{times}
\usepackage{graphicx}
\usepackage{amssymb}

\newtheorem{definition}{Definition}

\newtheorem{lem}{Lemma}[section]
\newtheorem{thm}{Theorem}[section]

\newtheorem{cor}{Corollary}[section]
\newtheorem{obs}{Observation}

\newcommand{\eop}{\hfill $\sqcap\!\!\!\!\sqcup$} 
\begin{document}
\global\def\refname{{\normalsize \it References:}}
\baselineskip 12.5pt
%
%
%
\title{\LARGE \bf The 3-rainbow index of graph operations}

\date{}

\author{\hspace*{-10pt}
\begin{minipage}[t]{2.7in} \normalsize \baselineskip 12.5pt
\centerline{TINGTING LIU}
\centerline{Tianjin University}
\centerline{Department of Mathematics}
\centerline{300072 Tianjin}
\centerline{CHINA}
\centerline{ttliu@tju.edu.cn}
\end{minipage} \kern 0in
\begin{minipage}[t]{2.7in} \normalsize \baselineskip 12.5pt
\centerline{YUMEI HU\footnote{corresponding author}}
\centerline{Tianjin University}
\centerline{Department of Mathematics}
\centerline{300072 Tianjin}
\centerline{CHINA}
\centerline{huyumei@tju.edu.cn}
\end{minipage}
\\ \\ \hspace*{-10pt}
\begin{minipage}[b]{6.9in} \normalsize
\baselineskip 12.5pt {\it Abstract:}
 A  tree $T$, in  an edge-colored graph $G$, is
called {\em a rainbow tree}  if no two edges of $T$ are assigned the
same color. A  {\em $k$-rainbow coloring }of $G$ is an edge coloring
of $G$ having the property that for every set $S$ of $k$ vertices of
$G$, there exists a rainbow tree $T$ in $G$ such that $S\subseteq
V(T)$. The minimum number of colors needed in a $k$-rainbow coloring
of $G$ is the {\em $k$-rainbow  index of $G$}, denoted by
$rx_k(G)$.  Graph operations, both binary and unary, are an interesting subject, which can be used to understand structures of graphs. In this paper, we will study the $3$-rainbow index with  respect to three important graph product operations (namely cartesian product, strong product, lexicographic product) and other graph operations. In this direction, we firstly show if $G^*=G_1\Box G_2\cdots\Box G_k$~($k\geq 2$), where each $G_i$ is connected, then $rx_3(G^*)\leq \sum _{i=1}^{k} rx_3(G_i)$. Moreover, we also present a condition and show the above equality holds if every graph $G_i~(1\leq i\leq k)$ meets the condition. As a corollary, we obtain
an upper bound for the 3-rainbow index of strong product. Secondly, we discuss the 3-rainbow index of the lexicographic graph $G[H]$  for connected graphs $G$ and $H$. The proofs are constructive and hence yield the sharp bound. Finally, we
consider the relationship between the 3-rainbow index of original graphs and other simple graph operations : the join of $G$ and $H$, split a vertex of  a graph and subdivide an edge.
\\ [4mm] {\it Key--Words:}
 $3$-rainbow index; cartesian product; strong product; lexicographic product.
\end{minipage}
\vspace{-10pt}}
\maketitle

\thispagestyle{empty} \pagestyle{empty}
\section{Introduction}
\label{S1} \vspace{-4pt}
All graphs considered in this paper are simple, connected and
undirected. We follow the terminology and notation of Bondy and
Murty \cite{bondy2008graph}. Let $G$ be a nontrivial connected graph
of order $n$ on which is defined an edge coloring, where adjacent
edges may be the same color. A path $P$ is a {\em rainbow path} if
no two edges of $P$ are colored the same. The graph $G$ is {\em
rainbow connected }if $G$ contains a $u$-$v$ rainbow path for every
pair $u,v$ of distinct vertices of $G$. If by coloring $c$ the graph $G$ is
rainbow connected , the coloring $c$ is called a rainbow coloring of $G$.
The {\em rainbow connection
number} $rc(G)$ of $G$, introduced by Chartrand et al. in
\cite{Chartrand2008graph}, is the minimum number of colors that
results in a rainbow connected graph $G$.

Rainbow connection has an interesting application for the secure transfer of
classified information between agencies (cf.~\cite{Ericksen2007}). Although the information needs to be protected since it is vital to national security, procedures must be in place that permit access between appropriate parties. This two fold issues can be addressed by assigning information transfer paths between agencies  which may have other agencies as intermediaries while requiring a large enough number of passwords and firewalls  that is prohibitive to intruders, yet small enough to manage (that is, enough so that one or more paths between every pair of agencies have no password repeated). An immediate question arises: What is the minimum number of passwords or firewalls needed that allows one or more secure paths between every two agencies so
that the passwords along each path are distinct? This situation can be modeled by a graph and studied by the means of rainbow coloring.

Later, another generalization of rainbow connection number was introduced by Chartrand
et al.\cite{Chartrand2009graph} in 2009. A tree $T$ is a {\em rainbow tree} if no two edges of $T$ are
colored the same. Let $k$ be a fixed integer with $2\leq k\leq n$.
An edge coloring of $G$ is called a {\em $k$-rainbow coloring} if
for every set $S$ of $k$ vertices of $G$, there exists a rainbow
tree in $G$ containing the vertices of $S$. The {\em $k$-rainbow
index} $rx_k(G)$ of $G$ is the minimum number of colors needed in a
$k$-rainbow coloring of $G$. It is obvious that $rc(G)=rx_2(G)$.
A tree $T$  is called a {\em concise} tree if $T$ contains $S$ and $T-v$ is not  a tree containing $S$, where $v$ is any vertex of $T$. In the paper, we suppose the tree containing $S$ be concise.
Since if the given tree $T$ is not concise, we can get a concise tree by deleting some vertices from $T$.

As we know, the diameter is a natural lower bound of the rainbow
connection number. Similarly, we consider the Steiner diameter in
this paper, which is a nice generalization of the concept of
diameter. The {\em Steiner distance} $d(S)$ of  a set $S$ of
vertices in $G$ is the minimum size of a tree in $G$ containing $S$.
Such a tree is called a {\em Steiner S-tree} or simply a {\em
Steiner tree}. The {\em $k$-Steiner diameter $sdiam_k(G)$} of $G$ is
the maximum Steiner distance of $S$ among all sets $S$ with $k$
vertices in $G$. The $k$-Steiner diameter  provides a lower
bound for the  $k$-rainbow index of $G$, i.e., $sdiam_k(G)\leq
rx_k(G)$. It follows, for every nontrivial connected graph $G$ of
order $n$, that
$$
rx_2(G)\leq rx_3(G)\leq \cdots \leq rx_k(G).
$$

For general $k$,
Chartrand et al. \cite{Chartrand2009graph} determined the $k$-rainbow index of
trees and  cycles. They obtained the following theorems.

\begin{thm}\cite{Chartrand2009graph}\label{thm3}
Let $T$ be a tree of order  $n\geq 3$. For each integer $k$ with $3\leq k
\leq n$, $$ rx_k(T)=n-1.$$
\end{thm}

\begin{thm}\cite{Chartrand2009graph}\label{thm5}
For  integers $k$ and $n$ with $3\leq k\leq n$,
 \[
      rx_k(C_n)=\left\{
      \begin {array}{lll}
      n-2, &\mbox{ if~ $k=3$ and $n\geq 4$;}\\
      n-1, &\mbox{ if~ $k=n=3$ or $4\leq k\leq n$.}\\
      \end {array}
      \right.
\]
\end {thm}

In the paper, we focus our attention on $rx_3(G)$.
For 3-rainbow index of a graph, Chartrand et al. \cite{Chartrand2009graph}
derive the exact value for the complete graphs.

\begin{thm}\cite{Chartrand2009graph}\label{thm4}
For any integer $n\geq 3$,
 \[
      rx_3(K_n)=\left\{
      \begin {array}{lll}
      2, &\mbox{ if~ $3\leq n\leq 5$;}\\
      3, &\mbox{ if~ $n\geq 6$;}\\
      \end {array}
      \right.
\]
\end {thm}

Chakraborty  et al. \cite{chakraborty2009hardness} showed that
computing the rainbow connection number of a graph is NP-hard. So it
is also NP-hard to compute $k$-rainbow index of a connected graph.
 For rainbow connection number $rc(G)$, people aim to give nice upper bounds
for this parameter, especially sharp upper bounds, according
to some parameters of the graph $G$ \cite{Chandran2011dominating,Li2012connected,LiShiSun,caro2008rainbow}.

 Many researchers have paid more attention to  rainbow connection number of some graph products \cite{M.B,SK,T,LS, m-2}. There is one way to bound the rainbow connection number of a graph product by the rainbow connection number of the operand graphs.   Li and Sun \cite{m-2} adopted the method  to study rainbow connection number with respect to Cartesian product and lexicographic product. They  got the following conclusions.

\begin{thm}\cite{m-2}\label{thm0}
Let $G^*=G_1\Box G_2\cdots\Box G_k$~($k\geq 2$), where each $G_i$ is connected, then

 $$rc(G^*)\leq \sum _{i=1}^{k} rc(G_i)$$
Moreover, if $rc(G_i)=diam(G_i)$ for each $G_i$, then the equality holds.
\end{thm}

\begin{thm}\cite{m-2}\label{thm0}
If $G$ and $H$ are two graphs and $G$ is connected, then we have \\
1. if  $H$ is complete, then

$$rc(G[H])\leq rc(G).$$
In particular, if $diam(G)=rc(G)$, then $rc(G[H])=rc(G)$.\\
2. if  $H$ is not complete,then

$$rc(G[H])\leq rc(G)+1.$$
In particular, if $diam(G)=rc(G)$, then $diam$$(G[H])=2$ if $G$ is complete and $rc(G)\leq diam(G)+1$ if $G$ is not complete.
\end{thm}

In this paper, we  study the $3$-rainbow  index with respect to three important graph product operations (namely cartesian product, lexicographic product and strong product) and other operations of
graphs. Moreover, we present the class of graphs  which obtain the upper bounds.

\subsection{Preliminaries}

 We use $V(G)$,~$E(G)$ for the set of vertices and edges of $G$, respectively.
For any subset $X$ of $V(G)$, let $G[X]$ be the subgraph induced by $X$, and $E[X]$ the edge set of $G[X]$; Similarly, for any subset $E'$ of $E(G)$, let $G[E']$ be the subgraph induced
by $E'$. For any two disjoint subsets $X$,~$Y$ of $V(G)$, we use $G[X,Y]$ to denote
the bipartite graph with vertex  set $X\cup Y$ and edge set $E[X,Y]=\{uv\in E(G)|u\in X, v\in Y\}$. The {\em distance} between two vertices $u$ and $v$ in $G$ is the length of a shortest path between them and is denoted by $d_G(u,v)$. The distance between a vertex $u$ and  a path $P$ is the shortest distance between $u$ and the vertices in $P$. Given a graph $G$, the eccentricity of a vertex, $v \in V (G)$ is given by $ecc(v) = max\{d_G(v, u): u\in V (G)\}$. The diameter of $G$ is defined as $diam(G) = max\{ecc(v): v\in V (G)\}$. The length of a path is the number of edges in that path.
The length of a tree $T$ is the numbers of edges in that tree, denoted by $size(T)$.
$G\setminus e$  denotes the graph obtained by deleting an edge $e$ from the  graph $G$ but leaving the vertices and the remaining edges intact. $G-v$ denotes the graph obtained by deleting the vertex $v$ together with all the edges incident with $v$ in $G$.

\begin{definition}(The Cartesian Product)
Given two graphs $G$ and $H$, the  Cartesian product of $G$ and $H$, denoted by $G\Box H$,
is defined as follows: $V (G\Box H)=V (G) \times V (H)$. Two distinct vertices $(g_1, h_1)$ and $(g_2, h_2)$ of $G\Box H$ are adjacent if and
only if either $g_1=g_2$ and $h_1 h_2\in E(H)$ or $h_1 = h_2$ and $g_1g_2\in E(G)$.
\end{definition}

\begin{definition}(The Lexicographic Product)
The  Lexicographic Product $G[H]$ of graphs $G$ and $H$ has the vertex set   $V(G[H])=V (G)\times V(H)$. Two vertices $(g_1, h_1), (g_2, h_2)$ are adjacent if $g_1g_2 \in  E(G)$, or
if $g_1 = g_2$ and $h_1h_2\in E(H)$.
\end{definition}

\begin{definition}(The Strong Product)
The Strong Product $G\boxtimes H$ of  graphs $G$ and $H$ is the graph with $V (G \boxtimes H)=V (G)\times V(H)$.
Two distinct vertices $(g_1, h_1)$ and $(g_2, h_2)$ of $G\boxtimes H$ are adjacent whenever  $g_1 = g_2$ and $h_1h_2 \in E(H)$ or $h_1 = h_2$ and $g_1g_2 \in E(G)$ or
 $g_1g_2 \in E(G)$ and $h_1h_2 \in E(H)$.
\end{definition}

 Clearly,  the resultant
graph is isomorphic to $G$ (respectively $H$) if $H = K_1$ (respectively $G = K_1$). Therefore,  we suppose $V(G)\geq 2$ and $V(H)\geq 2$ when  studying the 3-rainbow index of  these three graph products.

\begin{definition}(The union of graphs)
The union of two graphs, by starting with a disjoint union of two graphs $G$ and $H$ and adding edges joining every vertex of $G$ to every vertex of  $H$, the resultant graph is the {\em join } of $G$ and $H$, denoted by $G\vee H$.
\end{definition}

\begin{definition}(To split a vertex)
To {\em split} a vertex  $v$ of a graph $G$ is to replace  $v$ by two adjacent vertices
$v_1$ and $v_2$, and to replace each  edge incident to $v$ by an edge incident to
either $v_1$ or $v_2$~(but not both), the other end of the edge remaining unchanged.
\end{definition}
\subsection{Some basic observations}
\vspace{-4pt}
It is easy to see that if the graph $H$ has a $3$-rainbow coloring with $rx_3(H)$ colors, then the graph $G$, which is obtained from $H$ by adding  some edges to $H$,  also  has a $3$-rainbow coloring with $rx_3(H)$ colors  since  the new edges of $G$ can be colored arbitrarily with the colors used in $H$. So we have:

\begin{obs}\label{obs1}
Let $G$ and $H$ be  connected graphs and $H$ be a spanning subgraph of $G$. Then
$rx_3(G)\leq rx_3(H)$.
\end{obs}
To verify a 3-rainbow index, we need to find a rainbow tree containing any set of  three vertices. So it is necessary to know the structure of concise trees.
Next  we consider the structure of concise trees $T$ containing three vertices, which will be very useful in the sequel.

\begin{obs}\label{obs2}
Let $G$ be a connected graph and $S=\{v_1,v_2,v_3\}\subseteq V(G)$.  If ~$T$ is  a concise tree containing $S$, then $T$ belongs to exactly one of Type $I$ and Type $II$( see Figure 1).\\
Type $I$:  $T$ is a path such that one vertex of $S$ as its origin, one of $S$ as its terminus, other vertex of $S$ as its internal vertex.\\
Type $II$: $T$ is a tree obtained from  the star $S_3$ by replacing each edge of $S_3$ with a path $P$.
\end{obs}
\begin{figure}[h,t,b,p]
\begin{center}
\includegraphics[width=8cm]{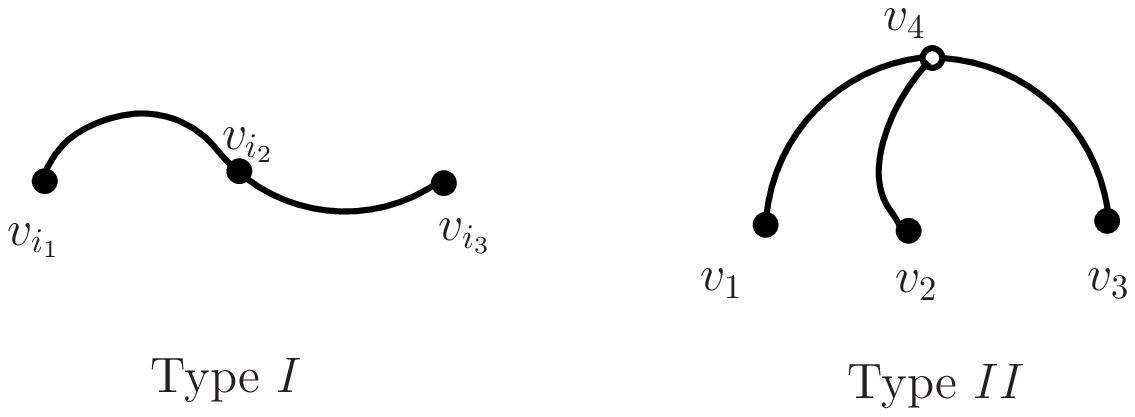}\\
Figure 1: Two types of concise trees, where $\{v_{i_1},v_{i_2},v_{i_3}\}=\{v_1,v_2,v_3\}$, $v_4\in V(G)$
\end{center}
\end{figure}
\noindent
{\bf Proof:} \
Firstly, we claim that the leaves of $T$ belong to $S$. Since if there exists a leaf $v$ such that $v\notin S$, then we can get the more minimal tree $T'=T-v$ containing $S$, a contradiction. Thus the $T$ has at most three leaves. If the $T$ has exactly two leaves, then it is easy to verify that $T$ is a path. In this case,  $T$ belongs to Type $I$. Otherwise there is a $v_1v_2$-path $P$ in $T$ such that $v_3\notin P$. Since $T$ is connected, there a path $P'$ in $T$ connecting $v_3$ and $P$.  Let $v_4$ be the vertex of $P'$ such that $d_T(v_3,v_4)$=$d_T(v_3,P)$. Then we get $T\supseteq P\cup P'$. On the other hand, we know,                                                                $P\cup P'$ is a tree containing $S$. Furthermore, since $T$ is a concise tree, $T=P\cup P'$, which belongs to Type $II$.\eop

\section{Cartesian product}
\vspace{-4pt}
In this section, we do some research on the relationship between the $3$-rainbow index of the original graphs and that of the cartesian products.
Recall that  the  Cartesian product of $G$ and $H$, denoted by $G\Box H$,
is defined as follows: $V (G\Box H)=V (G) \times V (H)$. Two distinct vertices $(g_1, h_1)$ and $(g_2, h_2)$ of $G\Box H$ are adjacent if and
only if either $g_1=g_2$ and $h_1 h_2\in E(H)$ or $h_1 = h_2$ and $g_1g_2\in E(G)$.
Let $V(G)=\{g_i\}_{i\in [s]}$, $V(H)=\{h_j\}_{j\in [t]}$.
Note that $H_i=G\Box H[\{(g_i,h_j)\}_{j\in [t]}]\cong H, G_j=G\Box H[\{(g_i,h_j)\}_{i\in[s]}]\cong G$. Any edge $(g_i,h_{j_1})(g_i, h_{j_2})$ of $H_i$ corresponds to edge $h_{j_1} h_{j_2}$ of $H$  and  $(g_{i_1},h_j)(g_{i_2}, h_j)$ of $G_j$ corresponds to
edge $g_{i_1} g_{i_2}$ of $G$.
For the sake of  our results, we give some useful and fundamental conclusions about the Cartesian product.

\begin{lem}\label{lem1}\cite{W}
The Cartesian product of two graphs is connected if and only if these two graphs
are both connected.
\end{lem}
\begin{lem}\label{lem2}\cite{W}
The Cartesian product is associative.
\end{lem}
\begin{lem}\label{lem3}\cite{W}
Let $(g_1,h_1)$ and $(g_2,h_2)$ be arbitrary vertices of  the Cartesian product $G\Box H$. Then $$d_{G\Box H}((g_1,h_1),(g_2,h_2))=d_G(g_1,g_2)+d_H(h_1,h_2).$$
\end{lem}
With the aid of Observation \ref{obs2} and above  Lemmas, we derive the following lemma, which is useful to show the sharpness of our main result.
\begin{lem}\label{lem4}
Let $G^*=G_1\Box G_2\cdots\Box G_k$~($k\geq 2$), where each $G_i$ is connected. Then $$Sdiam_3(G^*)=\sum _{i=1}^{k} Sdiam_3(G_i).$$
\end{lem}

\noindent
{\bf Proof:} \ We first prove the conclusion holds for the case $k=2$. Let $G=G_1$, $H=G_2$, $V(G)=\{g_i\}_{i\in [s]}$, $V(H)=\{h_j\}_{j\in [t]}$, $V(G^*)=\{g_i, h_j\}_{i\in [s],j\in [t]}=\{v_{i,j}\}_{i\in [s],j\in [t]}$. Let $S=\{(g_1, h_1),(g_2, h_2),(g_3, h_3)\}, S_1=\{g_1, g_2, g_3\}, S_2=\{h_1, h_2, h_3\}$ be a set of any three vertices of $V(G^*)$,~$V(G)$,~$V(H)$, respectively. Suppose that
$T$, $T_1$ and $T_2$  be   Steiner trees containing $S$, $S_1$, $S_2$, respectively. Next, we only need to show $size(T)$=$size(T_1)$+$size(T_2)$.

On the one hand, by the definition of the   Cartesian product of graphs, each edge of $G^*$ is exactly one element of $\{H_i, G_j\}$, $i\in [s], j\in [t]$. Then we can regard $T$ as the union $G'$ and $H'$, where $G'$ is induced by all the edges of $G_j\cap T$,~ $j\in [t]$, $H'$ is induced by all the  edges of $H_i\cap T$,~$i\in [s]$. Let $G''$ and $H''$ be the graphs induced by the corresponding edges of all edges of $G_j\cap T$ and $H_i\cap T$($i\in [s],j\in[t]$) in $G$ and $H$, respectively. Clearly, $G''$ and $H''$ are connected and containing $S_1$ and $S_2$, respectively. Hence, we have, $size(T)$=$size(G')$+$size(H')$=$size(G'')$+$size(H'')$ $\geq$ $size(T_1)$+$size(T_2)$.

On the other hand, we try to construct a tree $T'$ containing $S$ with $size(T')=$ $size(T_1)$+$size(T_2)$.  Notice that, for every subgraph in $G$ (or $H$), we can find the corresponding subgraph in any  copy $G_j$ ( or $H_i$).
If $T_1$ or  $T_2$  belongs to Type $I$, without loss of generality, say $T_1=P_1\cup P_2$, where $P_1$ is the path connecting $g_{i_1}$ and $g_{i_2}$, $P_2$ is the path connecting $g_{i_2}$ and $g_{i_3}$, $\{g_{i_1},g_{i_2},g_{i_3}\}=\{g_1,g_2,g_3\}$. We can find a tree $T'=P_1'\cup T_2'\cup P_2'$ containing $S$, where  the path $P_1'$ is the corresponding path of $P_1$ in $G_{i_1}$ and the path $P_2'$ is the corresponding path of $P_2$ in $G_{i_3}$, the tree $T_2'$ is the corresponding tree of $T_2$ in $H_{i_2}$, (see Figure 2).
\begin{figure}[h,t,b,p]
\begin{center}
\includegraphics[width=9cm]{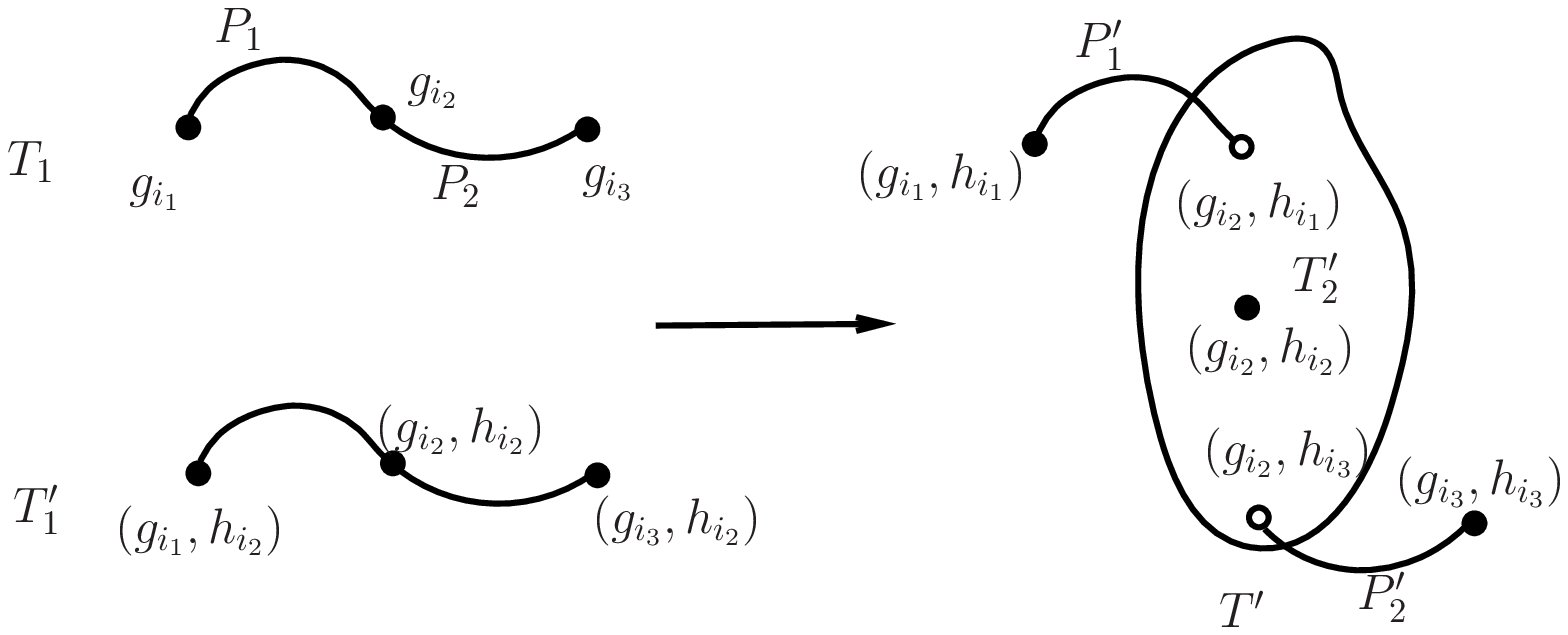}\\
Figure 2 : $T_1$ belongs to Type $I$
\end{center}
\end{figure}
If not, that is to say, $T_1$, $T_2$ belong to Type $II$, we suppose $T_1=P_1\cup P_2\cup P_3$, where $P_i$ is the path connecting $g_4$ and $g_i$ ($1\leq i\leq 3$), $g_4$ is other vertex of $G$ except the vertices of $S_1$.  Then the tree $T'=P_1'\cup P_2' \cup P_3'\cup T_2'$ containing $S$ can also  be found in $G\Box H$, where $P_i'$ is the corresponding path of $P_i$ in $G_i$~($1\leq i \leq 3$), the $T_2'$ is the corresponding tree of $T_2$ in $H_4$ (see Figure $3$). Thus, $size(T)\leq $ $size (T')$=$size(T_1)$+$size(T_2)$.
\begin{figure}[h,t,b,p]
\begin{center}
\includegraphics[width=8cm]{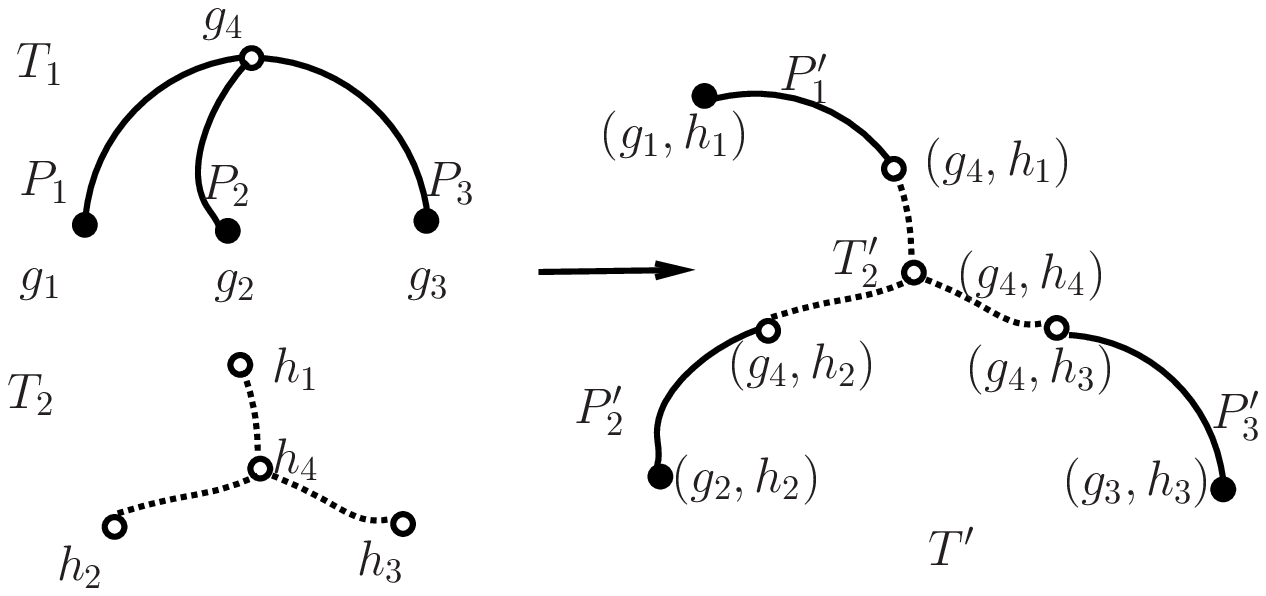}\\
Figure 3 : $T_1$ and $T_2$ belong to Type $II$.
\end{center}
\end{figure}

So we get  $size(T)$=$size(T_1)$+$size(T_2)$. Hence, $Sdiam_3(G_1\Box G_2)$= $Sdiam_3(G_1)$+$Sdiam_3(G_2)$. By Lemma \ref{lem2}, $Sdiam_3(G^*)$= $Sdiam_3(G_1\Box G_2$ $\Box \cdots $ $\Box G_{k-1})$+$Sdiam_3(G_k)$=$\sum _{i=1}^{k}$ $Sdiam_3(G_i)$.\eop

\begin{thm}\label{thm1}
Let $G^*=G_1\Box G_2\cdots\Box G_k$~($k\geq 2$), where each $G_i$ is connected, then

$$rx_3(G^*)\leq \sum _{i=1}^{k} rx_3(G_i)$$
Moreover, if $rx_3(G_i)=Sdiam_3(G_i)$ for each $G_i$, then the equality holds.
\end{thm}
\noindent
{\bf Proof:} \
We first show the conclusion holds for the case $k=2$. Let $G=G_1$, $H=G_2$, $V(G)=\{g_i\}_{i\in [s]}$, $V(H)=\{h_j\}_{j\in [t]}$, $V(G^*)=\{g_i, h_j\}_{i\in [s],j\in [t]}=\{v_{i,j}\}_{i\in [s],j\in [t]}$. Since $G$ and $H$ are connected,
$G^*$ is connected by Lemma \ref{lem1}. For example, Figure $4$ shows the case for
$G=P_4$ and $H=P_3$.

\begin{figure}[h,t,b,p]
\begin{center}
\includegraphics[width=6cm]{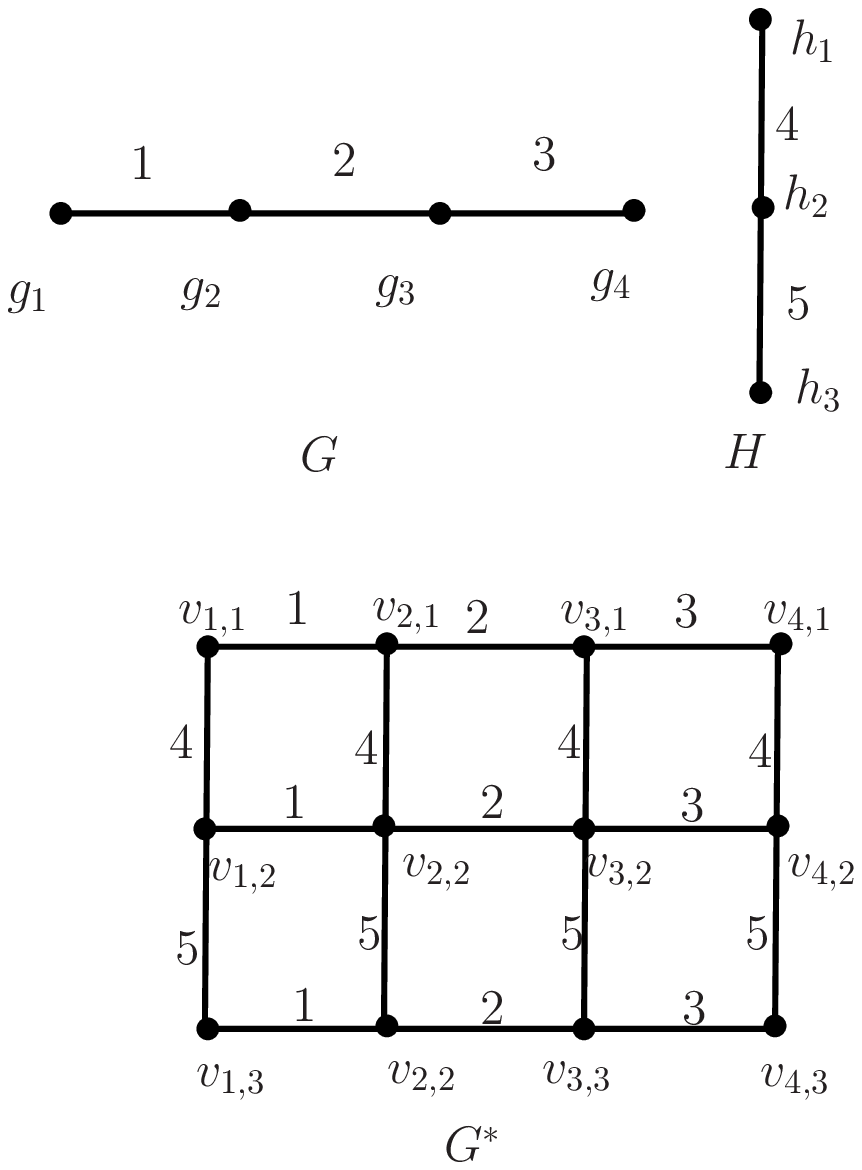}\\
Figure 4 : An example in Theorem \ref{thm1}.
\end{center}
\end{figure}

Since  for an edge $v_{i_1,j_1}v_{i_2,j_2}\in G^*$, we have $i_1=i_2$ or $j_1=j_2$;
if the former, then $v_{i_1,j_1}v_{i_1,j_2}\in H_{i_1}$, otherwise,
$v_{i_1,j_1}v_{i_2,j_1}\in G_{j_1}$. Hence, we only give a coloring of each graph $G_j~(j \in [t])$ and $H_i~(i\in [s])$.

We give $G$ a $3$-rainbow coloring with $rx_3(G)$ colors (see Figure 4 in which $G$ obtains a $3$-rainbow coloring with colors 1, 2, 3), and $H$ a $3$-rainbow coloring with $rx_3(H)$ fresh colors (see Figure 4 in which $H$ obtains a 3-rainbow coloring with other two fresh colors, 4, 5). Then we color edges of $G^*$ as follow:
if the edge belongs to some $H_i$, then assign the edge with the same color with its corresponding edge of $H$ (for example, edge $v_{1,1}v_{1,2}$ belong to $H_1$ and corresponds to the edge $h_1h_2$ in $H$, so it receives the color $4$), otherwise, the edge belongs to some $G_j$, then assign the edge with the same color with its corresponding edge of $G$. Now we will show that the given coloring is $3$-rainbow
coloring of $G^*$. It suffices to show that for every set $S$ of three vertices of $G^*$, there is a rainbow tree containing $S$. Let $S=\{(g_1,h_1),(g_2,h_2),(g_3,h_3)\}$. we distinguish three cases:

\textbf{Case 1} ~~ The vertices of $S$ lie in some $G_j$ (or $H_i$), where $i,j \in \{1,2,3\}$

That is, $g_1=g_2=g_3$ or $h_1=h_2=h_3$, without loss of generality,  we say, $g_1=g_2=g_3$. Under the given coloring of $H$, we can find a rainbow tree $T$ containing $ h_1,~h_2,~h_3$ in $H$. By the strategy of the above coloring, the corresponding tree $T'$  of $T$ in $H_1$ is also rainbow and contains $S$.

\textbf{Case 2} ~~ The vertices of $S$ lie in two different  copies $G_j'$ ,$G_j''$ (or $H_i'$, $H_i''$). where $j',~ j''\in \{1,~2,~3\}$(or $i',i''\in \{1~,2,~ 3\}$.

Without loss of generality, we assume $g_1=g_2\neq g_3$. Note that if a coloring is  $3$-rainbow coloring, then it is also rainbow coloring, that is, there is a rainbow path connecting any two vertices of graphs. If $h_1\neq h_2\neq h_3$~($h_1=h_3\neq h_2$ or $h_2=h_3\neq h_1$), we can find a rainbow tree  $T_1$ in $H$ containing $h_1,h_2,h_3$~($h_1,h_2$). By the strategy of coloring, we can find a rainbow tree $T_1'$ in $H_1$ containing $\{v_{1,1},v_{2,2},v_{1,3},\}$~($\{v_{1,1},v_{2,2}\}$).
So we can find a  rainbow path $P_1'$ in $G_3$ connecting $v_{1,3}$~($v_{1,1}$ or $v_{2,2}$) and $v_{3,3}$. Thus there is a rainbow tree $T=T_1'\cup P_1'$  in $G\Box H$ containing $S$.

\textbf{Case 3} ~~ The vertices of $S$ lie in three different copies  $G_1$,~$G_2$,~ $G_3$ and  $H_1$, ~$H_2$,~$H_3$.

 Let $T_1$ be a rainbow tree containing $g_1,g_2,g_3$ and $T_2$ be a
 rainbow tree containing $h_1,h_2,h_3$.

 If $T_1$ or $T_2$ belongs to
 Type $I$, say $T_1$, let $T_1=P_1\cup P_2$. Then the tree $T=P_1'\cup T_2'\cup P_2'$ containing $S$ can be constructed by  the way of Figure $2$. And by the character of the given coloring, the tree $T$ is a rainbow tree.

 If $T_1$ and $T_2$ belong to Type $II$, let $T_1=P_1\cup P_2\cup P_3$.  Then the tree $T=P_1'\cup P_2'\cup P_3'\cup T_2'$ can also be obtained by the way of Figure $3$. Furthermore, it is easy to see that the it is
 also a rainbow tree.

 Since we use $rx_3(G)+rx_3(H)$ colors totally, we have $rx_3(G^*)\leq rx_3(G)+rx_3(H)$.
 From Lemma \ref{lem4}, if $rx_3(G)=Sdiam_3(G)$~ and ~$rx_3(H)=Sdiam_3(H)$, then
 $Sdiam_3(G^*)=Sdiam_3(G)+Sdiam_3(H)=rx_3(G)+rx_3(H)\geq rx_3(G^*)$. On the other hand, $Sdiam_3(G^*)\leq rx_3(G^*)$, so the conclusion holds for $k=2$.

 For general $k$, by the Lemma \ref{lem2}, $rx_3(G^*)=rx_3(G_1\Box G_2\Box \cdots \Box G_{k-1}\Box G_k)\leq rx_3(G_1\Box G_2\Box $ $\cdots \Box G_{k-1})+rx_3(G_k)\leq
 \sum _{i=1}^{k}rx_3(G_i)$. Moreover, if $rx_3(G)=Sdiam_3(G_i)$ for each $G_i$, then
 $rx_3(G^*)\geq Sdiam_3(G^*)=\sum _{i=1}^{k} Sdiam_3(G_i)=\sum _{i=1}^{k} rx_3(G_i)\geq rx_3(G^*)$. So if $rx_3(G_i)=Sdiam_3(G_i)$ for each $G_i$, then the equality holds.\eop

\begin{cor} \label{cor1}
Let $G=P_{n_1}\Box P_{n_2}\Box \cdots \Box P_{n_k}$, where $P_{n_i}$ is a
path with $n_i$ vertices ($1\leq i \leq k$). Then

$$rx_3(G)=\sum_{i=1}^k n_i-k.$$
\end{cor}

\noindent
{\bf Proof:} \
 For  every path $P_{n_i}$, by Theorem \ref{thm3}, we have $Sdiam_3(P_{n_i})=rx_3(P_{n_i})=n_i-1$. Thus, by the Theorem \ref{thm1}, $rx_3(G)=\sum _{i=1}^{k} rx_3(P_{n_i})=\sum_{i=1}^k n_i-k$. \eop

Recall that the strong product $G \boxtimes H$ of graphs $G$ and $H$ has the vertex set
$V(G)\times V(H)$. Two vertices $(g_1, h_1)$ and $(g_2, h_2)$ are adjacent whenever $g_1=g_2 $  and  $h_1 h_2\in E(H)$  or  $h_1 = h_2$  and  $g_1g_2\in E(G)$ or $g_1g_2\in E(G)$  and  $h_1h_2\in E(H)$. By the definition, the graph $G\Box H$ is the spanning subgraph of the graph $G \boxtimes H$ for any graphs $G$ and $H$. With the help of Observation \ref{obs1}, then we have the following result.
\begin{cor}\label{cor2}
Let $\overline{G^*}$
=$G_1\boxtimes G_2\boxtimes\cdots\boxtimes G_k$,
~$(k\geq 2)$,
where each  $G_i~(1\leq i\leq k)$ is connected. Then we have

$$rx_3(\overline{G^*})\leq \sum_{i=1}^k rx_3(G_i).$$
\end{cor}
\section{Lexicographic Product}
\vspace{-4pt}
Recall that the lexicographic product $G[H]$ of graphs $G$ and $H$ has the vertex set $V(G[H]) = V(G)\times V(H)$. Two vertices $(g_1, h_1), (g_2
, h_2)$ are adjacent if $g_1g_2 \in  E(G)$, or
if $g_1= g_2$ and $h_1h_2\in E(H)$.  By definition, $G[H]$ can be obtained from
$G$ by submitting  a copy $H_1$ for every $g_1\in V(G)$ and by joining all vertices
of $H_1$ with all vertices of $H_2$ if $g_1g_2\in E(G)$.

In this section, we consider the relationship between 3-rainbow index of the original graphs and their lexicographic product.  Since the rainbow connection and 3-rainbow index is only defined in connected graphs, it is nature to assume the original graphs are connected.
Note that if $V(G)=1$ (or $V(H)=1$), then $G[H]$=$H$ (or $G$). So in the following discussion, we suppose $V(G)\geq 2$ and $V(H)\geq 2$.
By definition, if $G$ and $H$ are complete, then $G[H]$ is also complete.

So for some special cases of $G$ and $H$, we have the following lemma.
\begin{lem}\label{lem5}
If  $G,H\cong K_2$, then $$rx_3(G[H])=2.$$
If $G$ and $H$ are complete with $V(G)\geq 3$ or $V(H)\geq 3$, then $$rx_3(G[H])=3.$$
\end{lem}
\noindent
{\bf Proof:} \
If  $G, H\cong K_2$, then $G[H]$=$K_4$. Hence, we have $rx_3(G[H])=2$ by Theorem \ref{thm4}.  If $G$ and $H$ are complete with $V(G)\geq 3$ or $V(H)\geq 3$, then $G[H]=K_n$~($n\geq 6$). We get immediately $rx_3(G[H])=3$ from the Theorem \ref{thm4}.\eop

For the remaining cases, we obtain the following theorem.
\begin{thm}\label{thm2}
Let $G$ and $H$ be two connected graphs with $V(G)\geq 2$, $V(H)\geq 2$, and at least one of $G$, $H$ be not complete. Then\\
 $$rx_3(G[H])\leq rx_3(G)+rc(H).$$
In particular, if $diam(G)=rx_3(G)$, and $H$ is complete, then the equality holds.
\end{thm}

\noindent
{\bf Proof:} \
Let $V(G)=\{g_i\}_{i\in [s]}$, $V(H)=\{h_j\}_{j\in [t]}$, $V(G[H])=\{g_i, h_j\}_{i\in [s],j\in [t]}=\{v_{i,j}\}_{i\in [s],j\in [t]}$. Let $S=\{(g_1,h_1),(g_2,h_2),(g_3,h_3)\}$ be any three different vertices of $G[H]$.
 We derive the theorem from two parts:
1.~$V(H)=2$ and $G$ is not complete;
2.~$V(H)\geq 3$ and $G$ or $H$ is not complete.

1. If $V(H)=2$ and $G$ is not complete, we firstly give $G$ a $3$-rainbow coloring with $rx_3(G)$ colors. Then we can give $G[H]$ a $rx_3(G)$+1-edge coloring as follows: the edge belongs to
some $G_j$, then assign the edge with the same color with its corresponding edge
 in $G$. Otherwise, assign the edge a fresh color.

If $h_1=h_2=h_3$, then we can find a rainbow tree $T'$ containing $S$ , which is the corresponding tree of $T$  containing $g_1,~g_2,~g_3$ in $G_1$. Otherwise the vertices of $S$ lie in two different graphs $G_1$ and $G_2$. Without loss of generality, we suppose $h_1=h_3\neq h_2$.
In this case, $(g_1,h_1),(g_3,h_3)\in G_1$, $(g_2,h_2)\in G_2$.
 Then we can find the corresponding vertex  $(g_2,h_1)$~(or $(g_1,h_1)$ ~or~$(g_3,h_3)$) of $(g_2,h_2)$ in $H_1$ and a rainbow tree $T'$ containing $(g_1,h_1),(g_3,h_3)$ and  $(g_2,h_1)$~(or $\emptyset$). Clearly,
 there is a rainbow tree $T=T'\cup e$ containing $S$, where $e=(g_2,h_2) (g_2,h_1)$~(or $(g_1,h_1)$ or $(g_3,h_3)$). Hence the above coloring is $3$-rainbow coloring of $G[H]$. So $rx_3(G[H])\leq rx_3(G)+1=rx_3(G)+rc(H)$.

2. Let $c_1=\{0,1,\cdots,rx_3(G)-1\}$ be a $3$-rainbow coloring of $G$. Let $c_2$ be a rainbow coloring of $H$ using $rc(H)$ fresh colors. For every $h_j\in H$ color the copy $G_j$ the same as $G$.
By the same way, there is a rainbow tree containing any three vertices $(g_1,h_i),(g_2,h_i),(g_3,h_i)\in V(G[H])$. Every edge of the form $(g_1,h_1)(g_2,h_2)$ get color $k+1 $ mod($rx_3(G))$, where $g_1g_2\in E(G)$, $h_1\neq h_2$, and $c_1(g_1g_2)=k$. Finally, color edges from $H_i$ the same as $H$ such that any two vertices $(g_i,h_j)(g_i,h_k)$
 are connected by a rainbow path. The figure $5$ shows an example of the coloring.

\begin{figure}[h,t,b,p]
\begin{center}
\includegraphics[width=7cm]{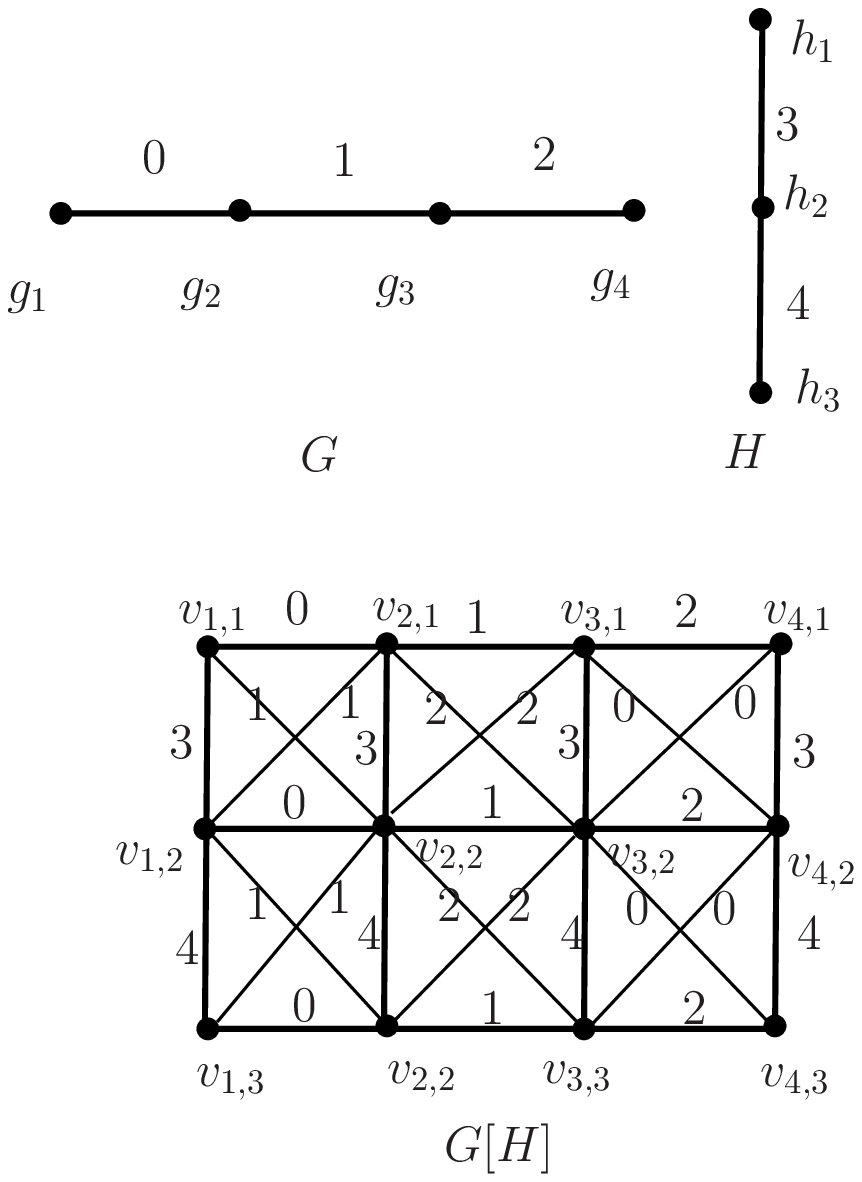}\\
Figure 5 : An example in Theorem \ref{thm2}. 2.
\end{center}
\end{figure}

Now we show the above coloring is $3$-rainbow coloring of $G[H]$.  We distinguish the following three cases.

\textbf{Case 1} ~~$g_1=g_2=g_3$

Since $G$ is a connected graph, there exists an edge $g_1g_4\in E(G), g_4\in V(G)$. Then
we can find a rainbow path $P$ connecting $(g_2,h_2)(g_1,h_1)$ in $H_1$, which uses the colors of $H$. By the coloring of strategy, the tree $T=P\cup v_{1,1}v_{4,1}\cup v_{4,1}v_{3,3}$ is a rainbow tree containing $S$.

\textbf{Case 2} ~~$g_1=g_2\neq g_3$ or  $g_1=g_3\neq g_2$ or  $g_2=g_3\neq g_1$

 Without loss of generality, we assume $g_1=g_2\neq g_3$.

 \textbf{Subcase 2.1}~~$h_1=h_3$ (or $h_2=h_3$)

  Then $T=P_1\cup P_2$ is a rainbow tree containing $S$, where $P_1$ is a rainbow path connecting $(g_1,h_1)$ and $(g_2,h_2)$ in $H_1$, $P_2$ is a rainbow path connecting $(g_1,h_1)$ (or $(g_2,h_2)$) and $(g_3,h_3)$ in $G_3$.

  \textbf{Subcases 2.2}~~ $h_1\neq h_2\neq h_3$

  As we know, there is a rainbow path $P_1$  connecting $g_3$ and $g_1$ in $G$. The case that  $P_1$=$g_3g_1$ is trivial, so we assume  $P_1$=$g_3g_1'$$g_2',$$\cdots,g_k'g_1$, $g_i'\in V(G)$~$(1\leq i \leq k)$. We claim that  $P_1'=(g_3,h_3)(g_1', h_2)(g_2',h_3)(g_3',h_2),\cdots, (g_k',u)(g_1,h_1)$ is a rainbow path connecting $(g_3,h_3)$ and $(g_1,h_1)$, where $u=h_3$ if $k$ is even and $u=h_2$ otherwise.
  It is easy to see that the path only use the edge of the form $(g_i,h_j)(g_j,h_l)$, where $g_ig_j\in E(G)$, $h_j\neq h_l$. By the character of
  coloring, the path is also  a rainbow path and only uses the colors of $G$.
  Thus,  there is  a rainbow tree $T=P'\cup P_2$ containing $S$, where  $P_2$ is a rainbow path connecting $(g_1,h_1)$ and $(g_2,h_2)$
 in $H_1$.

 \textbf{Case 3}~~ $g_1\neq g_2\neq g_3$

 \textbf{Subcase 3.1}~~$h_1=h_2=h_3$

  Then the $S$ lie in the  copy $G_1$. So by the given coloring, we can claim there is a rainbow tree $T$ containing $S$.

 \textbf{Subcase 3.2}~~$h_1=h_2\neq h_3$ or $h_1=h_3\neq h_2,$ or $ h_2=h_3\neq h_1$

 We suppose $h_1=h_2\neq h_3$. In this case, we first find the corresponding vertex
  $(g_3,h_1)$ of $(g_3,h_3)$ in $G_1$. Then  there is a rainbow tree  $T'$ containing
  $(g_1,h_1)(g_2,h_2)(g_3,h_1)$ in $G_1$ and a rainbow path $P$ connecting $(g_3,h_1)(g_3,h_3)$ in $H_3$. Thus, the rainbow tree $T=T'\cup P$ is our desire tree.

 \textbf{Subcase 3.3}~~$h_1\neq h_2\neq h_3$

   Suppose $T_1$ be a rainbow tree containing $g_1,g_2,g_3$.

   If $T_1$ or $T_2$ belongs to Type $I$, without loss of generality, we say $T_1$.
   In order to describe graphs simply,  we might suppose the leaves of $T_1$ are $g_1$ and $g_3$, $T_1=P_1\cup P_2$, where $P_1$ is a rainbow path connecting $g_1$ and $g_2$, $P_2$ is a rainbow path connecting $g_2$ and $g_3$. If $P_1$ or $P_2$ is an edge, it is trivial. So we  suppose $P_1=g_1g_1'g_2'$ $\cdots g_k'$ $g_2$ and $P_2=g_2g_1''g_2''$ $\cdots g_l''g_3$. Thus
    we can construct a rainbow tree $T_1'=P_1'\cup P_2'$ containing $S$, where
    $P_1'=(g_1,h_1)$$(g_1',h_3)(g_2',h_1)$
  $\cdots$ $(g_k',u)(g_2,h_2)$, $P_2'=(g_2,h_2)(g_1'', h_1)(g_2'',h_2)$
  $\cdots(g_l'',v)(g_3,h_3)$, $u=h_3$, if $k$ is odd, $u=h_1$ otherwise;
  $v=h_1$, if $l$ is odd; $v=h_2$ otherwise.

  If $T_1$ and $T_2$ belong to Type $II$, suppose $T_1=P_1\cup P_2\cup P_3$ and $T_2=Q_1\cup Q_2\cup Q_3$, where $P_i, Q_i~(1\leq i\leq 3)$ is a rainbow path connecting $g_4$ and $g_i$, $h_4$ and $h_i$. If $P_i~(1\leq i\leq 3)$ is an edge, then it is trivial. Now we suppose $P_i$ ($1\leq i\leq 3$) are
  not edges, then $P_1$=$g_4l_1'$$l_2'\cdots l_k'g_1$,~$P_2=g_4l_1''l_2''\cdots l_p'' g_2 $,~ $P_3$= $g_4l_1'''l_2'''\cdots l_q'''g_3$.
  Similarly, the corresponding rainbow tree $T_1'=P_1'\cup P_2'\cup P_3'$ can be obtained containing $S$, where $P_1'=(g_4,h_4)(l_1',h_2)$
 $(l_2',h_4)\cdots(l_k',u_1)$
 $(g_1,h_1)$,~ $P_2'=(g_4,h_4)$$(l_1'',h_3)$
 $(l_2'',h_4)\cdots(l_p'',u_2)$
 $(g_2,h_2)$, $P_3= (g_4,h_4)(l_1''',h_2)$
 $(l_2''',h_4)\cdots(l_q''',u_3)(g_3,h_3)$, $u_1,u_3=h_2,$
 $u_2=h_3$ if $k,p,q$ is odd, $u_1,u_2,u_3=h_4$, otherwise.

 From the above discussion, we have, the given coloring is $3$-rainbow coloring
 and we use $rx_3(G)+rc(H)$ colors totally. Thus, $rx_3(G[H])\leq rx_3(G)+rc(H)$.

 If $diam(G)=rx_3(G)$, and $H$ is complete, then $rx_3(G[H])\leq rx_3(G)+rc(H)=diam(G)+1$. On the other hand, let $g,g'\in V(G)$ such that
 $d_G(g,g')=diam(G)$. Let $S=\{(g',h),(g,h)(g,h')\}$. By the Lemma \ref{lem3},it is easy to check that the tree containing $S$ has
  size  at least $diam(G)+1$. So $rx_3(G[H])\geq Sdiam_3(G[H])\geq diam(G)+1$.
Thus,  $rx_3(G[H])=rx_3(G)+rc(H)$.\eop

\section{Other graph operations}
\vspace{-4pt}
We first consider the union of two graphs. Recall that the union of two graphs, by starting with a disjoint union of two graphs $G$ and $H$ and adding edges jointing every vertex of $G$ to every vertex of  $H$, the resultant graph is the  join  of $G$ and $H$, denoted by $G\vee H$.  Note that if $E(G)=\emptyset$ and $E(H)=\emptyset$, then the resultant graph is complete bipartite graph. So we need
 some results about the 3-rainbow index of complete bipartite graph. Li et al. got the following theorem for regular complete bipartite graphs $K_{r,r}$.

\begin{lem}\label{lem6}\cite{Li2013}
For integer $r$ with $r\geq3$, $rx_3(K_{r,r}) = 3$.
\end{lem}

For complete bipartite graph, we obtained the following Lemmas.
\begin{lem}\label{lem7}\cite{LH1}
For any  complete bipartite graphs $K_{s,t}$ with $3\leq s \leq t$,
$rx_3(K_{s,t})\leq min \{6,s+t-3\}$, and the bound is tight.
\end{lem}

In the proof of above Lemma \ref{lem7}, we showed the claim that
 for any $s\geq 3$, $t\geq 2\times 6^s$, $rx_3(K_{s,t})=6$.
\begin{lem}\label{lem8}\cite{LH2}
For any integer $t\geq 2$,
 \[
      rx_3(K_{2,t})=\left\{
      \begin {array}{lll}

      2, &\mbox{ if ~~$t=2$;}\\
      3, &\mbox{ if ~~$t=3,4$;}\\
      4, &\mbox{ if ~~$5\leq t\leq8$;}\\
      5, &\mbox{ if ~~$9\leq t\leq 20$;}\\
      k, &\mbox{ if ~~$C_{k-1}^2+1\leq t\leq C_k^2$,~($k\geq 6$).}\\
      \end {array}
      \right.
\]
\end{lem}

Then, we derive the relationship between the $3$-rainbow index of the original two
graphs and that of their join graph. Note that if $G$ and $H$ are both complete graphs, then $G\vee H$ is also the complete graph.  By the Theorem \ref{thm4}, $rx_3(G\vee H)=3$ if $|V(G)|$+$|V(H)|$$\geq 6$; $rx_3(G\vee H)=2$ if $|V(G)|$+$|V(H)|\leq 5$. So we consider the remaining cases in following theorem.

\begin{thm}\label{thm6}
If $G$, $H$ are connected and at least one of $G$, $H$ are not complete, with $|V(G)|=s$, $|V(H)|=t$,~$s\leq t$, then we have\\
1. if ~~$s=1$, then $$rx_3(G\vee H)\leq rx_3(H)+1.  $$
2. if ~~$2=s\leq t$, then $$rx_3(G\vee H)\leq min\{rc(H)+3, rx_3(K_{2,t})\}.$$
3. if ~~$3\leq s \leq t$, then $$rx_3(G\vee H)\leq min\{c_1+1, rx_3(K_{s,t})\}$$
Where $c_1=max\{rx_3(G),rx_3(H)\}$.\\
In particular, if $s=t\geq 3$, then $rx_3(G\vee H)=rx_3(K_{s,t})=3$.
\end{thm}
\noindent
{\bf Proof:} \
Let $G'=G\vee H$, $V(G')=V_1\cup V_2$ such that $G'[V_1]\cong G$, $G'[V_2]\cong H$,
where $V_1=\{v_1,v_2,\cdots, v_s\}$, $V_2=\{u_1,u_2,\cdots,u_t\}$.

1. If $s=1$, then  $G'[V_1]$ is singleton vertex, we give an edge coloring of $G'$
as follows : we first give a 3-rainbow coloring of $G'[V_2]$ using $rx_3(H)$ colors. And for the other edges, that is, elements of $E[V_1,V_2]$, we use a fresh color.
It is easy to show the above coloring of $G'$ is 3-rainbow coloring.

2. If $2=s\leq t$, then $G'[V_1,V_2]\cong K_{2,t}$ is a spanning subgraph of $G'$.
We have $rx_3(G')\leq rx_3(G'[V_1,V_2])=rx_3(K_{2,t})$. On the other hand, we
give an edge coloring of $G'$
as follows: we first color the edges of the subgraph $G'[V_2]$  with $rc(H)$
colors such that it is rainbow connected; we give the elements of $E[V_1,V_2]$  incident with $v_i$($i=1,2$) with color $rc(H)+i$ ($i=1,2$); for the element of $G'[V_1]$, we use  a fresh color $rc(H)+3$.
It is easy to show the above coloring of $G'$ is 3-rainbow coloring. Thus, we have
$rx_3(G\vee H)\leq min\{rc(H)+3, rx_3(K_{2,t})\}$.

3. If ~~$3\leq s \leq t$, by Lemma \ref{lem7}, we have $rx_3(G')\leq rx_3(G'[V_1,V_2])=rx_3(K_{s,t})$, similarly.
On the other hand, we color the edges of $G'$ as follows: we first color
the edges of the subgraph $G'[V_i]$ with $c_1$
 colors such that it is 3-rainbow coloring of $G'[V_i]$~($i=1,2$).
 For the rest edges, that is,  elements of $E[V_1,V_2]$, we use a fresh color $c_1+1$.
 It is  easy to verify that the coloring is a $3$-rainbow coloring. Thus, we
 get $rx_3(G\vee H)\leq min\{rx_3(K_{s,t}),c_1+1\}$.

 If $s=t\geq 3$, by Lemma \ref{lem6}, then $rx_3(G')\leq rx_3(K_{s,s})=3$;
 On the other hand, by Observation \ref{obs1} and Theorem \ref{thm4}, $rx_3(G')\geq rx_3(K_{s+t})=3$, so the conclusion holds.

 Note that $rx_3(K_{2,t})$ may be larger than $rc(H)+3$; for example, we choose $H\cong K_t\setminus e~(t\geq 21)$. Then $rx_3(K_{2,t})>5=rc(H)+3$ by Lemma \ref{lem8}. But $rx_3(K_{2,t})$ is not always larger than
 $rc(H)+3$; for example, we choose $H\cong P_t~$, then $rx_3(K_{2,t})<t+2=rc(H)+3$.
 Moreover, $rx_3(K_{s,t})~(3\leq s <t)$ may be larger than $max\{rx_3(G),rx_3(H)\}+1$, since we suppose $G\cong K_s \setminus e~(s\geq 3)$ and $H\cong K_t$, where $t\geq 2\times 6^s$. Then $rx_3(K_{s,t})=6 > max\{rx_3(G),rx_3(H)\}+1$. But $rx_3(K_{s,t})$ is not always larger than
 $max\{rx_3(G),rx_3(H)\}+1$. Similarly, for example, $G,H\cong P_s$~($s\geq 7$), we
 can get $max\{rx_3(G),rx_3(H)\}+1=s>6\geq rx_3(K_{s,t})$. So the bounds we give in the theorem are reasonable.\eop

 Recall that  to split $v$ of a graph $G$  is to replace  $v$ by two
 adjacent vertices $v_1$ and $v_2$ by an edge incident to either $v_1$
 or $v_2$ (but not both), the other end of the edge remaining unchanged. The Figure $6$ shows the operation of $G$. Let $N_G(v)$ be the neighbor sets of $v$. The set
is partitioned into two disjoint sets $N_1$ and $N_2$ such that
 $N_1$ and $N_2$ are the neighbor sets of $v_1$ and $v_2$ in the resultant
 graph, respectively.

\begin{figure}[h,t,b,p]
\begin{center}
\includegraphics[width=6cm]{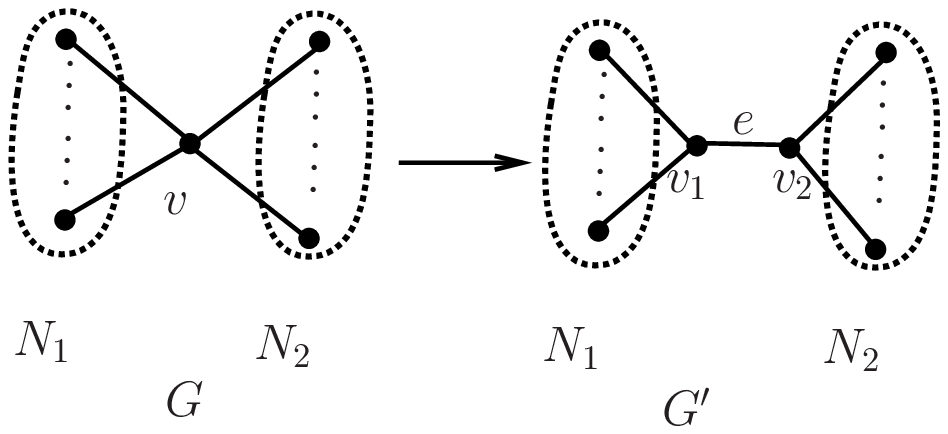}\\
Figure 6 : The operation  for vertex spliting.
\end{center}
\end{figure}
\begin{thm}\label{thm7}
 If $G$ is a connected graph and  $G'$ is obtained from $G$ by splitting a vertex $v$, then $$rx_3(G')\leq rx_3(G)+1.$$
\end{thm}
\noindent
{\bf Proof:} \
We first give $G$ a $3$-rainbow coloring with $rx_3(G)$ colors, then we give $G'$
a $rx_3(G)$+1-edge coloring as follows: we  give the edge $e=v_1v_2$ a color $rx_3(G)$+1; for any edge $uv_1\in G'$ with $uv_1\neq e$, let the color of $uv_1$
 be the same as that of $uv$ in $G$; for any edge $v_2w\in G'$ with
  $v_2w\neq e$, let the color of $v_2w$ be the same as that of $vw$ in $G$; color
  of the rest edges of $G'$ are the same as in $G$.
  Next, we will show the given coloring of $G'$ is a 3-rainbow coloring. It suffices to show that there is a rainbow tree containing any three vertices of $G'$. Let
  $S=\{x,y,z\}$.

  \textbf{Case 1}~~ Two  vertices of $S$ belongs to $\{v_1,v_2\}$, say $x=v_1$, $y=v_2$.

  By the above coloring, there a rainbow $v-z$ path $P:v=u_1,\cdots,u_t=z$. If $u_2\in N_1$, then $P':v_1,u_2,u_3,\cdots, u_t=z$ is a rainbow connecting $z$
   and $x(v_1)$. Thus, $T=P'\cup e$ is the rainbow tree containing $S$.
  If $u_2\in N_2$, it is similar to verify that there is a rainbow tree containing $S$.

  \textbf{Case 2}~~ Exactly one of $S$ belongs to $\{v_1,v_2\}$, say $x=v_1$.

   We know that, in graph $G$, there is a rainbow tree $T_1$ containing $y,z,v$.

   \textbf{subcase 2.1} ~~$d_{T_1}(v)=1$.

    Then there is an edge $uv\in E(T_1)$. If $u\in N_1$, the  tree  obtained from $T_1$ by replacing $v$ with $v_1$ is rainbow and contains $S$. If $u\in N_2$, the tree  obtained from $T_1$ by replacing $v$ with $v_2,~v_1$ is a rainbow
    tree containing $S$.

    \textbf{subcase 2.2}~~$d_{T_1}(v)\neq 1$.

     From the Observation \ref{obs2}, we claim $d_{T_1}(v)=2$. Let $u_1$ and $u_2$
     be  the two neighbors of $v$ in $T_1$. If $u_1$ and $u_2$ belong to the $N_1$, then let $T$ be obtained from $T_1$ by replacing  $v$  with $v_1$. If $u_1$ and $u_2$ belong to the $N_2$, then  we can find a rainbow tree $T=T_2\cup e$, where $T_2$ is obtained from $T_1$ by replacing $v$ with $v_2$. If $u_1$ and $u_2$ belong to the different $N_i~(i=1,2)$, then $T$ obtained from $T_1$ by
     replacing $v$ with subgraph $v_1v_2$ is rainbow.

   \textbf{Case 3}~~ None of vertices in $S$ belongs to $\{v_1,v_2\}$.

   We know that there is  a rainbow $T_3$ containing $S$ in $G$.
   If $v$ does not belong to $T_3$, then $T_3$ is also a rainbow tree containing $S$ in $G'$.

   If $v$ belong to the tree $T_3$, by the Observation \ref{obs2}, then $d_{T_3}(v)=2,~3$. Similar to the Subcase 2.2, we can find a rainbow tree
   containing $S$.

   So $G'$ receives a 3-rainbow coloring. Since we use  $rx_3(G)+1$ colors totally,
   then $rx_3(G')\leq rx_3(G)+1$.\eop

A special case of vertex splitting occurs when exactly one link is assigned to
 either $v_1$ or $v_2$. The resulting graph can be viewed as having
been obtained by subdividing an edge of the original graph, where to {\em subdivide} an edge is to delete $e$, add a new vertex $x$, and join $x$ to the ends of $e$. So
by Theorem \ref{thm7}, we have
\begin{cor}\label{cor3}
If $G$ is a connected graph, and $G'$ is obtained from $G$ by subdividing an edge
$e$, then $$rx_3(G')\leq rx_3(G)+1.$$
\end{cor}

\section{Conclusion}
\vspace{-4pt}
Rainbow connection number $rc(G)$~($rx_2(G)$) comes from the communication of information between agencies of government. $3$-rainbow index, $rx_3(G)$,  is a generalization of rainbow connection number. Chakraborty et al. have proved that computing $rc(G)$($rx_2(G)$) is NP-hard. Hence, To get the exact value for 3-rainbow index of general graph $G$ is also NP-hard. Thus, researchers tend to get some better upper for 3-rainbow index of some classes of graphs. Graph operations, both binary and unary, are  interesting subjects, which can be used to understand structures of graphs. In this paper, we will study the $3$-rainbow index with  respect to three important graph product operations (namely cartesian product, strong product, lexicographic product) and other graph operations. In this direction, we firstly show if $G^*=G_1\Box G_2\cdots\Box G_k$~($k\geq 2$), where each $G_i$ is connected, then $rx_3(G^*)\leq \sum _{i=1}^{k} rx_3(G_i)$. Moreover, we also present a condition and show the above equality holds if every graph $G_i~(1\leq i\leq k)$ meets the condition. As a corollary, we obtain
an upper bound for the 3-rainbow index of strong product. Secondly, we discuss the 3-rainbow index of the lexicographic graph $G[H]$  for connected graph $G$ and $H$. The proofs are constructive and hence yield the sharp bound. Finally, we
consider the relationship between the 3-rainbow index of original graphs and other simple graph operations : the join of $G$ and $H$, split a vertex of  a graph and subdivide an edge and get the upper bounds.

\vspace{10pt} \noindent
{\bf Acknowledgements:} \ The corresponding author,
Yumei Hu, is supported by NSFC No. 11001196.

\end{document}